\documentclass[12pt]{article}
\usepackage{amssymb}

\newcommand{\Z}{{\mathbb Z}}
\newcommand{\F}{{\mathbb F}}
\newcommand{\R}{{\mathbb R}}
\newcommand{\SUM}{\raisebox{-0.4ex}{\mbox{\Large $\Sigma$}}}
\newcommand{\eps}{\varepsilon}

\newtheorem{theorem}{Theorem}
\newtheorem{lemma}{Lemma}

\title{Arithmetic structures in smooth subsets of $\F_p$}

\author{Ernie Croot}

\begin{document}

\maketitle

\begin{abstract}  Suppose $S \subseteq \Z_N := \Z/N\Z$, where $N$ is a prime.
A well-studied question for various types of sets $S$, is that of whether for 
a particular sequence of integers $a_1,...,a_d$ satisfying $a_1 + \cdots + a_d = 0$,
$S$ contains solutions to the congruence
$$
a_1 s_1 + \cdots + a_d s_d\ \equiv\ 0\pmod{N}.
$$
In particular, the question for $d=3$ and $a_1 = a_2 = 1$ and $a_3 = -2$ is that 
of whether or not $S$ contains three-term arithmetic progressions.

In the present paper, rather than working with sets $S$ we work with functions 
$f : \Z_N \to [0,1]$, and consider the counting function
$$
\SUM_{x_1,...,x_d \in \Z_N \atop a_1 x_1 + \cdots + a_d x_d \equiv 0 \pmod{N}}
f(x_1)\cdots f(x_d).
$$
We show that if the function $f$ is ``sufficiently smooth'' --that is, 
the sum of squares of the ``small'' Fourier coefficients is sufficiently 
small -- then this counting function must be ``large''.  The proof is a generalization
of an arugment from an earlier paper of the author \cite{croot}, and it appears that
there are some close parallels with that proof and Green's ``arithmetic 
regularity lemma'' appearing in \cite{green} (see remark 3 at the end of the 
Introduction below).

One may think of this result as a statement about the circle method.  In that context
it says that, regardless of what the ``major arc'' contribution to the counting function
above happens to be, so long as along the ``minor arcs'' there are only very few places
where the exponential sum corresponding to $f$ (i.e. Fourier transform of $f$) can
be ``large'', then the counting function above much be ``large''.  So, one does not
even need to bother computing the contribution of the major arcs if one is only interested
in lower bounds for the counting function.  Of course to get asymptotics, the major
arcs would need to be evaluated precisely.

The theorem is proved in the following way.  First, in section 
\ref{separation_property_section} we precondition the Fourier coefficients of $f$
by applying a certain dilation function, so that the places $a$ where $|\hat f(a)|$
satisfy a certain technical conditional we call the ``separation property''.  
Next, we multiply $f$ by a certain ``smoothing function'', which will allow us to
transfer the problem of showing that our counting function is large, to an analogous
problem in $\Z_M$, where $M > N$ can be factored as $M = m_1 m_2$; if the
count in the analogous problem is ``large'', then so must be the counting problem in 
$\Z_N$.  Not just any number $M$ will do -- it must satisfy a certain property we
call the ``correspondance property'', and finding such $M$ will involve carefully
selecting $m_1$ and applying the ``separation property''.
Then, we apply ideas from \cite{croot} to replace our new smooth counting function
with one that is just as smooth, but also translation-invariant by the subgroup 
of $\Z_M$ consisting of multiples of $m_2$.  Finally, we show that the new 
counting function is ``large'', which means the same is true of all previous ones.
\end{abstract}

\section{Introduction}

For a prime $N$ we use the abbreviation $\Z_N\ :=\ \Z/N\Z$.  
Suppose that for $N \geq 2$ (we do not assume it is prime) we 
have a function
$$
f\ :\ \Z_N\ \to\ [0,1].
$$
In the present paper we will prove a theorem which says that if $f$ 
is ``sufficiently smooth'', in the sense that the sum of squares of the 
small Fourier  coefficients of $f$ are ``small'', then there are lots of arithmetic
structures on which $f$ is positive; for example, there will be lots of 
progression triples $n,n+d,n+2d$ where
$$
f(n) f(n+d) f(n+2d)\ >\ 0.
$$
In the context of the Circle Method, our theorem says the following:  For a certain 
class of additive problems, if one can show that along the ``minor arcs'' the 
measure of the places where the corresponding exponential sum is 
``not-too-small'', is itself ``small'', then there is no need to bother 
working out the contribution of the 
``major arcs'' (unless one wants asymptotic estimates), because one 
can show that regarless of what it is, one must have a large positive count
in the end.  

We will apply our theorem to prove, for example, that there are lots 
of three-term arithmetic progressions among certain sumsets and among 
the pseudoprimes (pseudoprimes
of the type considered by Goldston, Pintz and Yilidrim).   Although these results 
are already known, our result will give the same conclusion for quite large classes
of sets similar to pseudoprimes.

Before we state our main theorem, we need a few definitions.  First, for an $a \in \Z_N$
we define the Fourier transform
$$
\hat f(a)\ :=\ \SUM_n f(n) \omega^{an},\ {\rm where\ } \omega\ =\ e^{2\pi i /N}.
$$
We let $\lambda_1,...,\lambda_N$ be the Fourier coefficients of $f$ ordered so 
that 
$$
|\hat f(0)|\ =\ |\lambda_1|\ \geq\ |\lambda_2|\ \geq\ \cdots\ \geq\ |\lambda_N|.
$$
We let 
$$
\sigma^2\ :=\ \SUM_j |\lambda_j|^2\ =\ N \SUM_n |f(n)|^2\ \leq\ \theta N^2.
$$
This second equality follows from Parseval.

Our theorem is as follows.

\begin{theorem} \label{main_theorem}  Fix 
$$
d\ \geq\ 3,\ {\rm and\ } a_1,\ a_2,\ ...,\ a_d\ \in\ \Z \setminus \{0\},
$$
satisfying
$$
a_1 + \cdots + a_d\ =\ 0,
$$
and fix 
$$
0\ <\ \eps\ <\ 1.
$$
Then the following holds for all primes $N$ sufficiently large: 
Suppose 
$$
f\ :\ \Z_N\ \to\ [0,1]
$$ 
satisfies
$$
\SUM_n f(n)\ =\ \theta N\ >\ 0,
$$
and has the property that for some integer $k$ satisfying
$$
1000^{d \eps^{-1}} \theta^{-(\eps d)^{-1}}\log N\ \leq\ k\ \leq\ N^{1/11}, 
$$
we have
$$
\SUM_{k \leq j \leq N} |\lambda_j|^2\ <\ k^{-(4+10\eps)(d-2)} 
\sigma^2\ \leq\ k^{-(4 + 10\eps)(d-2)} \theta N^2;
$$

\noindent Then,
$$
\SUM_{x_1,...,x_d \in \Z_N \atop a_1 x_1 + \cdots + a_d x_d \equiv 0 \pmod{N}} 
f(x_1)\cdots f(x_d)\ >\ 10^{-1} 4^{-d} k^{-2(d-2)-2\eps(d-3/2)}\theta (\theta N)^{d-1}
$$
\end{theorem}

\noindent {\bf Remark 1.}  The upper bound we demand for $k$ can
be substantially improved, though it would take quite a bit of work to get it
above $N^{1/2}$, assuming this is even possible.
\bigskip

\noindent {\bf Remark 2.}  The lower bound we prove here has the general shape of what
we should expect:  Assuming that $f$ is the indicator function for some set,
there are $(\theta N)^{d-1}$ choices for $x_1,...,x_{d-1}$ such that 
$f(x_1)\cdots f(x_{d-1}) > 0$, and then we expect that 
only $\theta$ fraction of the values 
$$
x_d\ \equiv\ -a_d^{-1}(a_1 x_1 + \cdots + a_{d-1} x_{d-1}) \pmod{N}
$$
land in the set as well.  So, a {\it reasonable} lower bound should be 
$\theta (\theta N)^{d-1}$ for the final sum in the statement of the theorem.
\bigskip

\noindent {\bf Remark 3.}  There appears to be a relationship between
Theorem \ref{main_theorem} and ``triangle deletion'' ideas of 
Szmeredi and Ruzsa \cite{ruzsa} and of Green \cite{green}.  Indeed, one of the 
central ideas in the present paper appears to use some of the same types of
ingredients as those of Green's theorems from \cite{green} (which I discovered after
finally skimming Green's paper!).
In order to be more specific, it is worth looking at Green's ``arithmetic regularity
lemma'' in the case of $G := (\Z/2\Z)^n$.  In that proof he constructs a sequence of
smaller and smaller subgroups $H$ (which nonethelss are still quite large) 
of $G$, until one fairly large one is found 
that has certain nice ``regularity'' properties; and, these subgroups appear only to 
be definable in an iterative manner.  By having some precondition on the sum of
squares of the small Fourier coefficients of our basic starting function $f$, as we 
do in Theorem \ref{main_theorem} above (though imagine $f$ is defined on $G$), 
we can bypass this iterative process,
and can, in fact, just pick our $H$ randomly, and of fairly high dimension 
(a positive probability of the $H$ we could pick will work).  Furthermore,
this seems to work even when the functions that one uses have very low density
(that is, $\SUM_n f(n)$ is ``small'').  (I should
say that there only {\it appears} to be a connection between the two papers, 
as I have not thought about it in depth.)  There is still the 
problem of how to make the idea work 
modulo $N$.  If one tries to use Bohr neighborhoods as in Green's paper, one will
have lots of new technical complications to deal with; however, in our 
proof of Theorem \ref{main_theorem} we bypass
these problems by passing to another group $\Z_M$ that has a large additive 
subgroup with certain usable properties.
\bigskip

We now devote a new section to give two common examples of functions $f$ where the
sum is positive, at least in the case corresponding to three-term arithmetic progressions,
which is $k=3$ and $a_1 = a_2 = 1$ and $a_3 = -2$:  
The first example is sumsets, and the second is pseudoprimes.
For both of these examples one can establish the existence of such three-term arithmetic
progressions by other methods, so these examples are only meant to be suggestive of 
what types of results one can obtain from our theorem.

\section{Some types of sets where Theorem \ref{main_theorem} applies}

\subsection{Sumsets and arithmetic progressions}

Suppose that $S$ is a subset of $\Z_N$ having at least $N^{0.999}$ elements, which means
that its density is $\theta \sim N^{-0.001}$.  We will show
how Theorem \ref{main_theorem} implies that the $6$-fold sumset $S+S+S+S+S+S$ contains
a three-term arithmetic progression.  Note, however, that we know that just the 
$2$-fold sum $S+S$ contains three-term progressions by an elementary argument; so,
Theorem \ref{main_theorem} does not give anything new when we apply it to repeated
sumsets.  

To see how to prove this fact about $6$-fold sumsets using Theorem \ref{main_theorem},
first define
$$
f(n)\ :=\ |S|^{-5} (S*S*S*S*S*S)(n),
$$
which is supported exactly on the sumset $S+S+S+S+S+S$, and has size at most 
$1$ there.  If we let $\gamma_1,...,\gamma_N$ denote the Fourier coefficients of 
$S$, ordered so that 
$$
|S|\ =\ \gamma_1\ \geq\ |\gamma_2|\ \geq\ \cdots\ \geq\ |\gamma_N|,
$$
then we have that the Fourier coefficients of $f$ are $\lambda_1,...,\lambda_N$,
where
$$
\lambda_i\ =\ |S|^{-5} \gamma_i^6.
$$

We now give an upper bound for $|\gamma_k|$ by observing from Parseval that
$$
k |\gamma_k|^2\ \leq\ \SUM_i |\gamma_i|^2\ =\ N |S|.
$$

This then implies that
\begin{eqnarray*}
\SUM_{i \geq k} |\lambda_i|^2\ &\leq&\ |S|^{-10} 
|\gamma_k|^{10} \SUM_{i \geq 1} |\gamma_i|^2 \\
&\leq&\ k^{-5} |S|^{-4} N^6.
\end{eqnarray*}
Now for $\eps = 1/11$ we will have that
$$
\SUM_{i \geq k} |\lambda_i|^2\ <\ k^{-4-10\eps} (k^{-\eps} |S|^{-4}N^6). 
$$
So, for $k \sim \theta^{-66}$, which is certainly smaller than $N^{1/11}$, 
we will have that
$$
\SUM_{i \geq k} |\lambda_i|^2\ <\ k^{-4-10\eps} \SUM_i |\lambda_i|^2;
$$
and so, Theorem \ref{main_theorem} implies that $S$ contains a three-term arithmetic
progression on letting $d=3$ and $a_1 = a_2 = 1$ and $a_3 = -2$.  

\subsection{Pseudoprimes}
 
We define a pseudoprime in the sense of Goldston, Pintz and Yildirim \cite{goldston} as 
used in the work of Green and Tao \cite{green_tao}:  First, let $\delta > 0$ and then 
for all $n \leq N/2$, say, let
$$
f(n)\ :=\ {1 \over (\log N)^2 \max_{n \leq N/2} \tau(n)^2} 
\left ( \SUM_{d | n \atop d \leq N^\delta} \mu(d) \log(N/d) \right )^2.
$$
We note that if $n$ is a prime number, then 
$$
f(n)\ =\ {1 \over (\log N) \max_{n \leq N/2} \tau(n)^2}\ >\ N^{-o(1)},
$$
and regardless of whether $n$ is prime or not, we have that 
$$
f\ :\ \{ n \leq N/2 \}\ \to\ [0,1].
$$
We furthermore have that if we think of $f$ as a function on $\Z_N$ (in the obvious way),
then
$$
\hat f(0)\ \geq\ N^{1-o(1)},
$$
since there are $N^{1-o(1)}$ primes $\leq N/2$.  

What is not immediately obvious, but true, is that $f$ can be easily perturbed 
so that $\hat f$ is ``smooth enough'' for Theorem \ref{main_theorem} to imply
that there are lots of arithmetic structures where $f$ is positive, at least if our 
``truncation level'' $N^\delta$ in the definition of $f$ above is small enough.  
In order to see this, we first define 
$$
g(n)\ :=\ \SUM_{d | n \atop d \leq N^\delta} \mu(d) \log(N/d),
$$
and observe that $f(n)$ is the square of $g(n)$ up to a scalar factor of size 
$N^{-o(1)}$.  The Fourier transform of $g$ is given by
\begin{eqnarray*}
\hat g(a)\ &=&\ \SUM_{n \leq N/2} e^{2\pi i an/N} \SUM_{d | n \atop d \leq N^\delta} \mu(d) 
\log(N/d) \\
&=&\ \SUM_{d \leq N^\delta} \mu(d) \log(N/d) \SUM_{n \leq N/2 \atop d | n} e^{2\pi i an/N}.
\end{eqnarray*}
This inner sum at the end is a geometric series, and is a fairly 
``smooth'' function.  Unfortunately, it is not quite smooth enough for the 
particular way that we apply Theorem \ref{main_theorem}; so, we will need to
multiply it by a certain weighting function $w(n)$ to make it even smoother.  
The following standard, well-known technical lemma does this for us.

\begin{lemma}  For all $1 \leq d \leq N^{0.0001}$, 
there exists a weighting function $w_d : \Z_N \to [0,1]$ such that 
\bigskip

$\bullet$  The function $w_d$ is supported at most on the set of all integer multiples 
of $d$ lying in $[0,N/2]$ modulo $N$;
\bigskip

$\bullet$ all but at most $N^{0.005}$ points $a \in \Z_N$ we have that 
$$
\left | \SUM_{dn \leq N/2} w_d(dn) e^{2\pi i adn/N}  \right |\ <\ N^{-2};
$$
\bigskip

$\bullet$ and, $w_d(dn) = 1$ for $1000 d N^{0.999} < dn < N/2 - 1000 d N^{0.999}$.
\end{lemma}

\noindent {\bf Proof of the Lemma.}  In the proof we will just drop the subscript $d$
on $w_d$.  

First, let 
$$
X\ :=\ \lfloor N^{0.999} \rfloor,
$$
and then we define $h : \Z_N \to \R_{\geq 0}$ via its Fourier transform 
$$
\hat h(a)\ :=\ \left ( \SUM_{0 \leq n < X} e^{2\pi i ad n/N} \right )^{1000}.
$$
Then, we define $w$, also through its Fourier transform, as
$$
\hat w(a)\ =\ X^{-1000} \hat h(a) \SUM_{0 \leq j \leq N/2000dX} e^{2\pi i a j (1000 dX)/N}. 
$$

It is not difficult to see that  
$$
w\ :\ \Z_N\ \to\ [0,1]
$$
that 
$$
{\rm support}(w)\ \subseteq\ \{ d n\ :\ 0 \leq n \leq N/2d\},
$$
and that for 
$$
1000dX\ <\ dn\ <\ N/2 - 1000d X
$$
we have that
$$
w(dn)\ =\ 1.
$$
Furthermore, it is easy to give non-trivial upper bounds on the number of places 
$a$ where $\hat w$ is ``large'':  First, observe that we have the trivial upper bound
$$
|\hat w(a)|\ <\ X^{-1000} N |\hat h(a)|.
$$
So, the only places $a$ where $|\hat w(a)|$ could exceed $N^{-2}$ are those
where 
$$
|\hat h(a)|\ >\ X^{1000} N^{-3}.
$$
In other words,
$$
\left | \SUM_{0 \leq n < X} e^{2\pi i ad n/N} \right |\ >\ X N^{-3/1000}\ >\ (1/2)N^{0.996}.
$$
It is a simple matter to prove that there can be at most $N^{0.005}$ places $a$
having this property.

$\hfill$ $\blacksquare$
\bigskip

From this lemma we see that if we replace $g$ with $g_2$, where
\begin{equation} \label{g*}
\hat g_2(n)\ =\  \SUM_{d \leq N^\delta} \mu(d) \log(N/d) 
\SUM_{1 \leq dn \leq N/2} w_d(dn) e^{2\pi i adn/N},
\end{equation}
then by the third bullet in the lemma above we find that 
$$
g_2(m)\ =\ g(m)\ \ {\rm whenever\ \ } 1000 d N^{0.999}\ <\ m\ <\ N/2 - 1000 d N^{0.999}.
$$

We also have that $g_2$ can have only very few places $a$ where 
$|\hat g_2(a)| > N^{-1}$:  By the second bullet in the lemma above, along
with the definition (\ref{g*}), we see that the number of places $a$ where
$\hat g_2(a)$ can exceed $N^{-1}$ is at most the number of places $a$
where any one of the inner sums (for any $d \leq N^\delta$) of (\ref{g*})
exceeds $N^{-2}$.  The number of such places $a$ is clearly bounded from
above by
$$
N^{\delta + 0.005}.
$$

Next we show that the function 
$$
f_2(n)\ :=\ {g_2(n)^2 \over (\log N)^2 \max_{n \leq N/2} \tau(n)^2},
$$
also has only very few places $a$ where its Fourier transform is not too small: 
First, observe that
$$
\hat f_2(a)\ =\ {(g_2*g_2)(a) \over N  (\log N)^2 \max_{n \leq N/2} \tau(n)^2}\ =\ 
N^{-1-o(1)} \SUM_{x+y \equiv a \pmod{N}} \hat g_2(x) \hat g_2(y).
$$
In order for $a$ to be such that 
$$
|\hat f_2(a)|\ >\ N^{-1},
$$
we must have that $a \equiv x + y \pmod{N}$ where $x$ and $y$ are 
places where $|\hat g_2(x)|$ and $|\hat g_2(y)|$ exceed $N^{-1}$.  The number
of such $a$, then, is at most the square of the number of places $x$ where
$|\hat g_2(x)| > N^{-1}$.  It follows then that there are at most 
$$
N^{2 \delta + 0.01}
$$
such places $a$.  

We now pass to one more function $f_3(n)$ by performing yet one more 
level of smoothing.  The reason for this is that we only have that 
\begin{equation} \label{f2}
f_2(m)\ =\ f(m)\ \ {\rm whenever\ \ } 1000 N^{0.999+\delta}\ <\ m\ <\ N/2 - 1000 
N^{0.999+\delta}.
\end{equation}
Thus, we want to zero out the function $f_2(m)$ for values of $m$ that 
are close to $0$ or $N/2$, while still maintaining the fact that our function
has few not-so-small Fourier coefficients.  The function $f_3$ we will use is
given by
$$
f_3(m)\ :=\ f(m) w^*(m),
$$
where $w^*$ is to be defined via its Fourier transform as follows
$$
\hat w^*(a)\ =\ C e^{2\pi i a \lceil N/5\rceil /N} 
\left ( \SUM_{0 \leq n \leq N/5000} e^{2\pi i am/N} \right )^{1000},
$$  
where $C$ is chosen so that $\sup_m w^*(m) = 1$.  

Since $w^*$ is supported at most on $[N/5, 2N/5]$, which is well within the range
(\ref{f2}), we deduce that 
$$
f_3(m)\ >\ 0\ \ \Longrightarrow\ \ f(m)\ >\ 0.
$$
Furthermore, it is a routine calculation to show that 
$$
\hat f_3(0)\ =\ N^{1-o(1)},
$$
and that for all but at most 
$$
N^{2\delta + 0.02}
$$
places $a \in \Z_N$ we have that 
$$
|\hat f_3(a)|\ <\ N^{-1/2}.
$$
This clearly implies that if we let $\lambda_1,...,\lambda_N$ be
the Fourier coefficients of $f_3$, ordered so that 
$$
|\lambda_1|\ \geq\ |\lambda_2|\ \geq\ \cdots\ \geq\ |\lambda_N|,
$$
then for $k \sim N^{2 \delta + 0.02} < N^{1/11}$ (for small $\delta >0$) 
we will have that 
$$
\SUM_{j \geq k} |\lambda_j|^2\ \leq\ k^{-5} \SUM_i |\lambda_i|^2.
$$
Theorem \ref{main_theorem} then implies that (for 
$d = 3$ and $a_1 = a_2 = 1$ and $a_3 = -2$) there are lots of 
three-term progressions $m,m+t,m+2t$ such that 
$$
f_3(m) f_3(m+t) f_3(m+2t)\ >\ 0
$$
and therefore lots of progressions where
$$
f(m)f(m+t)f(m+2t)\ >\ 0.
$$

\section{Proof of Theorem \ref{main_theorem}}

Throughout the proof we will make use of the two parameters
$$
D\ :=\ 4 d \max_{1 \leq i \leq d} |a_i|,\ {\rm and\ } L\ :=\ \lfloor \log N \rfloor + 1.
$$

To prove the theorem we will move the problem from $\Z_N$ to an additive 
group $\Z_M$ having certain subgroups with useful properties.   We now work 
this out in the following subsection.

\subsection{Moving to another group}

We fix a prime number $m_2$ in advance that satisfies
\begin{equation} \label{m2_interval}
k^{2+2\eps}\ \leq\ m_2\ \leq\ 2 k^{2+2\eps},
\end{equation}
and then we will later find an integer $m_1$ satisfying
\begin{equation} \label{m1_interval}
k^{-2-\eps} N\ \leq\ m_1\ \leq\ 2 k^{-2-\eps} N,
\end{equation}
such that we can transfer our counting from from $\Z_N$ to $\Z_M$, where
$$
M\ =\ m_1 m_2\ \in\ [k^\eps N,\ 4 k^\eps N].
$$  
The advantage of making this transfer is that $\Z_M$ has a 
relatively small index subgroup (index $m_2$) consisting of the 
multiples of $m_1$ which we will later exploit.
In order to make this go through smoothly, however, we will need to not only
select $m_1$ very carefully, but will need to precondition the Fourier coefficients of
$f$, to get them to satisfy what we call the ``separation property''.

\subsubsection{Separation property of Fourier coefficients} 
\label{separation_property_section}

Let 
$$
\{b_1,...,b_k\}\ \subseteq\ \Z_N
$$
be the places satisfying
$$
\hat f(b_i)\ =\ \lambda_i.
$$

By replacing $f$ with the function
$$
f^*(n)\ :=\ f(qn),
$$
for an appropriate $1 \leq q \leq N-1$, we will show that we may assume that the Fourier
coefficients $b_1,...,b_k$ satisfy the following ``separation property'':
\bigskip

\noindent {\bf Separation property.}  All but a fraction 
$k^{-2}$ of the integers $m_1$ satisfying (\ref{m1_interval})
have the property that for every quadruple
$$
(b_i, b_j,a_u,a_v)\ \in\ \{b_1,...,b_k\}^2 \times \{a_1,...,a_d\}^2,
$$
we have that 
\begin{eqnarray} \label{separation_eitheror}
{\rm either\ }a_u b_i - a_v b_j\ \equiv\ 0 \pmod{N},\ {\rm or\ } 
\left | \left | {m_1(a_u b_i - a_v b_j) \over N} \right | \right |\ >\ {k^{4\eps} \over m_2}. 
\nonumber \\
\end{eqnarray}

\subsubsection{Proof that the separation property can be satisfied}

First, let us note that replacing $f$ with $f^*$ will not affect our weighted count
(weighted by $f$) of the number of solutions to
$$
a_1 x_1 + \cdots + a_d x_d \equiv 0 \pmod{N},
$$ 
since this congruence ``respects dilations'', in the sense that it holds if and only if
$$
a_1 (q x_1) + \cdots + a_d (q x_d) \equiv 0 \pmod{N},
$$
whenever $(q,N) = 1$.  That is to say, we will have that 
\begin{eqnarray*}
&& \SUM_{x_1,...,x_d \in \Z_N \atop a_1 x_1 + \cdots + a_d x_d \equiv 0 \pmod{N}}
f(x_1)\cdots f(x_d)\\
&&\hskip0.5in =\ \SUM_{x_1,...,x_d \in \Z_N \atop a_1 x_1 + \cdots + a_d x_d
\equiv 0 \pmod{M}} f^*(x_1)\cdots f^*(x_d).
\end{eqnarray*}
Also, we note that
$$
\hat f^*(qa)\ =\ \hat f(a),
$$
which means that in place of $\{b_1,...,b_k\}$ where $\hat f$ is of ``large size'',  
we can work with $\{qb_1,qb_2,...,qb_k\}$ where $\hat f^*$ has ``large size''.  
\bigskip

Our job now is to show that there exists a value for $q$ such that if we let 
$qb_1,...,qb_k$ stand in place of $b_1,...,b_k$, then the separation property
can be made to hold.  We will do this using some harmonic analysis, and
we begin by letting $J$ be the integers $m_1$ satisfying 
(\ref{m1_interval}), and letting
$$
K\ :=\ [-k^{4\eps} N/m_2,\ k^{4\eps} N/m_2] \cap \Z.
$$ 
By deleting at most one element we can make $K = -K$, and will assume this
is so.  Although $J$ and $K$ are defined as integer intervals, we will
think of them as subsets of $\Z_N$.

Next, suppose we fix a quadruple
\begin{equation} \label{quadruple1}
(b_i, b_j, a_u, a_v)\ \in\ \{b_1,...,b_k\}^2 \times \{a_1,...,a_d\}^2,
\end{equation}
satisfying
\begin{equation} \label{quadruple2}
a_u b_i - a_v b_j \not \equiv 0 \pmod{N}.
\end{equation}
Then, for a fixed $q$, the number of integers $m_1 \in J$ such that 
$$
m_1(a_u (qb_i) - a_v (qb_j))
$$ 
lies in the interval $K$ modulo $N$ is bounded from above by
\begin{equation} \label{C_def}
C\ =\ 3|K|^{-1} N^{-1} \SUM_a \hat J(ax) \hat K(-a)^2,
\end{equation}
where
$$
x\ \equiv\ a_u (qb_i) - a_v (qb_j) \pmod{N}.
$$
The reason for this is as follows:  First, if we write out
$$
\hat K(a)^2\ =\ \SUM_{-N/2 < n \leq N/2} \ell(n) e^{2\pi i an/N},
$$
then 
$$
\ell(n)\ =\ |K| - |n|,\ {\rm for\ } |n|\ \leq\ |K|.
$$
For $|K| \leq |n| \leq N/2$ the function $\ell(n)$ will have value $0$.  
Note that for $n \in K$ our formula for $\ell(n)$ implies that 
$$
\ell(n)\ >\ |K|/3,
$$
and this lower bound is the origin of the factor $3 |K|^{-1}$ appearing in our 
upper bound (\ref{C_def}) on our count for the number $m_1 \in J$ above.
\bigskip

Since $\hat J(b)$ and $\hat K(b)$ are geometric series, it is easy to prove that for
$$
-N/2\ <\ b\ \leq\ N/2,
$$
we have
\begin{equation} \label{geometric_bound}
|\hat J(b)|\ \leq\ \min( |J|,\ |N/b|),\ {\rm and\ } |\hat K(b)|\ \leq\ \min(|K|,\ |N/b|);
\end{equation}
in particular this means
$$
|\hat J(ax)|\ \leq\ \min(|J|,\ ||ax/N||^{-1}).
$$
An estimate we will need in a minute is 
$$
|K|^{-1}\SUM_a |\hat K(a)|^2\ <\ 20N,
$$
which can be proved by using the upper bound $|K|$ on the size of $|\hat K(a)|$
for when 
$$
|a|\ <\ 3 |K|^{-1} N,
$$
and then applying (\ref{geometric_bound}) for when 
$$
3|K|^{-1} N\ \leq\ |a|\ \leq\ N/2.
$$
\bigskip

Now, let us suppose for the time being that $x$ is such that 
\begin{equation} \label{goodx}
||ax/N||\ >\ N |J|^{-1}|K|^{-1},\ {\rm for\ all\ } 0\ <\ |a|\ \leq\  N^2 |K|^{-2}.
\end{equation}
Then we have that 
$$
C\ =\ 3N^{-1} |J| \cdot |K|\ +\ E_1\ +\ E_2,
$$
where 
\begin{eqnarray*}
|E_1|\ &\leq&\  3|K|^{-1} N^{-1} \sup_{|a| \leq N^2 |K|^{-2}  \atop a \neq 0} |\hat J(ax)| 
\SUM_a |\hat K(a)|^2\\
&<&\ 3N^{-1} |J| \cdot |K|.
\end{eqnarray*}
and where
\begin{eqnarray*}
|E_2|\ &\leq&\ 3 |K|^{-1} |J| N^{-1} \sum_{N^2 |K|^{-2} < |a| \leq N/2}
|\hat K(a)|^2  \\
&<&\ 3 |K|^{-1} |J| N \SUM_{|a| > N^2 |K|^{-2}} |a|^{-2} \\
&<&\ 6 N^{-1} |J| \cdot |K|.
\end{eqnarray*}
So, for $x$ satisfying (\ref{goodx}) we have that 
$$
C\ <\ 20N^{-1} |J| \cdot |K|.
$$
What this means is that all but a fraction $20N^{-1} |K| < k^{-2}$ 
of the integers in $m_1 \in J$ must satisfy  
$$
|| m_1 ( a_u (qb_i) - a_v (qb_j))/N | |\ >\ k^{4\eps}/ m_2.
$$
\bigskip

Recalling that 
$$
x\ \equiv\ a_u (qb_i) - a_v (q b_j) \pmod{N},
$$
we have that the number of values of $q$ that fail to satisfy 
the first inequality of (\ref{goodx}) for a particular non-zero 
value of $a$ is at most 
$$
2N^2 |J|^{-1} |K|^{-1}.
$$
So, the number $q$ failing to satisfy (\ref{goodx}) is at most 
\begin{eqnarray*}
(2N^2 |J|^{-1} |K|^{-1})(2N^2 |K|^{-2})\ &=&\ 
4 N^4 |J|^{-1} |K|^{-3} \\
&\leq&\ k^{2 - 3 \eps} m_2^3.
\end{eqnarray*}
So, the number of $q$ failing to satisfy (\ref{goodx}) for all quadruples 
satisfying (\ref{quadruple1}) and (\ref{quadruple2}) is, by (\ref{m2_interval}), 
at most
\begin{equation} \label{last_count}
(d^2 k^2) k^{2 - 3\eps} m_2^3\ \leq\ 8 d^2 k^{10 + 3\eps}. 
\end{equation}

Since $k < N^{1/11}$ and $0 < \eps < 1/3$, we have that the last quantity of 
(\ref{last_count}) is smaller than $N$, and therefore there exists $q$ such that 
we can make the separation property hold.
\bigskip

\subsubsection{The auxilliary function $g$}

Let
$$
I\ :=\ [-k^{-\eps} N,\ k^{-\eps} N] \cap \Z.
$$
Then, define, for a certain value of $u \in \Z_N$ to be decided in a moment, 
\begin{equation} \label{ell_formula}
\hat g(a)\ :=\ N^{-1} |I|^{-L + 1} \SUM_{-N/2 < b \leq N/2} e^{2\pi i ub/N} \hat f(b) 
\left ( \SUM_{n \in I} e^{2\pi i n (a/M - b/N)} \right )^L.
\end{equation}
We have that $\hat g(a)$ may be written as
$$
\hat g(a)\ :=\ \SUM_{-M/2 < n \leq M/2} f(n-u) w(n) e^{2\pi i a n/M},\ {\rm for\ some\ }
w\ :\ \Z \to [0,1],
$$
where here we are thinking of $f$ as a periodic mapping $f : \Z \to [0,1]$ 
having period $N$, instead of as a mapping $f : \Z_N \to [0,1]$.
Furthermore, the function $w$ satisfies
$$
\SUM_n w(n)\ =\ |I|.
$$

Note that by Fourier inversion this means that $g$ is supported at most on
the interval $(-N/2, N/2]$ modulo $M$, at least if $k^\eps > 2L$, and that 
$$
g(n)\ =\ f(n-u) w(n),\ {\rm for\ } -N/2\ <\ n\ \leq\ N/2.
$$
By simple averaging we have that there exists $u$ such that 
\begin{equation} \label{gmass}
\hat g(0)\ =\ \SUM_n f(n-u)w(n)\ \geq\ N^{-1} |I| \hat f(0)\ \geq\ k^{-\eps} \theta N,
\end{equation}
and we will assume that we have used any such $u$ to define our function $g$.

Since this function $g$ is only supported at most on the 
interval $(-N/2, N/2]$ modulo $M$, we have that if 
$$
-M/2\ <\ x_1,...,x_d\ \leq\ M/2
$$
satisfies
$$
a_1 x_1 + \cdots + a_d x_d\ \equiv\ 0 \pmod{M}\ \ {\rm and\ \ } 
g(x_1)\cdots g(x_d)\ >\ 0,
$$
then, in fact, $x_1,...,x_d$ must be confined to the smaller interval 
$(-N/2, N/2]$, which then implies that
$$
a_1 x_1 + \cdots + a_d x_d\ =\ 0\ \ {\rm in\ }\Z,
$$
and therefore this also holds modulo $N$.

So,
\begin{eqnarray} \label{fg}
&& \SUM_{x_1,...,x_d \in \Z_N \atop a_1 x_1 + \cdots + a_d x_d \equiv 0 \pmod{N}} 
f(x_1)\cdots f(x_d) \nonumber \\
&& \hskip0.5in\ \geq\ \SUM_{x_1,...,x_d \in \Z_M \atop a_1 x_1 + \cdots + a_d x_d \equiv 0 
\pmod{M}} g(x_1)\cdots g(x_d).
\end{eqnarray}
\bigskip

\subsubsection{The sizes of the Fourier coefficients of $g$}

Next we need to better understand the size of the Fourier coefficients of
$g$.

As this last factor of (\ref{ell_formula}) is a geometric series, we have that
\begin{eqnarray} \label{g_fourier}
|\hat g(a)|\ &\leq&\ N^{-1} |I|^{-L+1} \SUM_{-N/2 < b \leq N/2} |\hat f(b)|
\min \left ( |I|^L,\ 2^L |1 - e^{2\pi i (a/M - b/N)}|^{-L} \right ) \nonumber \\
&\leq&\ N^{-1} |I|^{-L+1} \SUM_{-N/2 < b \leq N/2} |\hat f(b)| \min(|I|^L,\ 
|\sin(\pi(a/M - b/N))|^{-L}) \nonumber \\
&\leq&\ N^{-1} |I|^{-L+1} \SUM_{-N/2 < b \leq N/2} |\hat f(b)| 
\min(|I|^L,\ 2^{-L} ||a/M - b/N||^{-L}). \nonumber \\
\end{eqnarray}

We now use this to get some handle on the places $a$ where $|\hat g(a)|$ is 
``large'':  First, let 
$$
S\ =\ \{b_1,...,b_k\},\ {\rm and\ then\ let\ } S^c\ :=\ \Z_N \setminus S.
$$
Then, we know that
$$
b \in S^c\ \ \Longrightarrow\ \ |\hat f(b)|\ \leq\ |\lambda_k|\ \leq\ 
k^{-(2 + 4 \eps)(d-2)} \theta^{1/2} N.
$$

Next let $X$ denote the set of all integers $a$ such that
$$
-M/2\ <\ a\ \leq\ M/2,
$$
having the property
\begin{equation} \label{xmeaning}
{\rm for\ all\ }i=1,...,k,\ ||a/M - b_i/N||\ >\ k^{3\eps}/M.
\end{equation}
Note that
$$
|X^c|\ <\ 3k^{1+3\eps},\ {\rm and\ }|X|\ >\ M - 3k^{1+3\eps}.
$$

For such $a\in X$ we will have that 
\begin{equation} \label{b_bad}
2^{-L}|I|^{-L+1} ||a/M - b_i/N||^{-L}\ <\ N^{-10}.
\end{equation}
for $N$ sufficiently large, since $k > L$.

Now, using (\ref{g_fourier}), consider the sum
\begin{eqnarray} \label{g_upper}
\SUM_{a \in X} |\hat g(a)|^2\ &\leq&\ 
N^{-2} |I|^{-2L+2} \SUM_{a \in X} \SUM_{-N/2 \leq c, c' \leq N/2}
 |\hat f(c)|\cdot | \hat f(c')| \nonumber \\
&&\hskip0.5in \times\ \min(|I|^L, \ 2^{-L} ||a/M - c/N||^{-L}) \nonumber \\
&&\hskip0.5in \times\ \min(|I|^L,\ 2^{-L} ||a/M - c'/N||^{-L}). 
\end{eqnarray}
In this sum, we first observe that the contribution of those pairs $(c,c')$ that
are ``far apart'' is very small.  Specifically, if 
\begin{equation} \label{dlc}
||c/N - c'/N||\ >\ k^{3\eps}/N,
\end{equation}
then we will have that for every $a \in \Z_M$,
$$
{\rm either\ } ||a/M - c/N||\ \ {\rm or\ \ } ||a/M - c'/N||\ >\ k^{3\eps}/2N.
$$
If either of these occurs, say the first one occurs, then we will have that
$$
\min(|I|,\ 2^{-L} |I|^{-L+1} ||a/M - c/N||^{-L})\ <\\ N^{-10}.
$$
So, the total contribution of all such terms to (\ref{g_upper}) will be at most 
$N^{-6}$.   

From this and (\ref{g_upper}) it is not difficult to see that this implies
\begin{eqnarray*}
\SUM_{a \in X} |\hat g(a)|^2\ &\leq&\ 
N^{-2} |I|^{-2L+2} \SUM_{a \in X} \SUM_{-N/2 \leq c,c' \leq N/2 \atop 
||c/N - c'/N|| \leq k^{3\eps}/N}
 |\hat f(c)|\cdot | \hat f(c')| \nonumber \\
&&\hskip0.25in \times\ \min(|I|^L, \ 2^{-L} ||a/M - c/N||^{-L}) \nonumber \\
&&\hskip0.25in \times\ \min(|I|^L,\ 2^{-L} ||a/M - c'/N||^{-L}) + N^{-6}.
\end{eqnarray*}

Now, if $c$ or $c'$ equals $b_1,...,b_k$, then from the fact that 
$a \in X$, one or the other
of these last two factors, when multiplied by $|I|^{-L+1}$, 
will be smaller than $N^{-10}$, making the total contribution
of those terms very small.  On the other hand, if $c$ and $c'$ both fail to equal 
$b_1,...,b_k$, then we get that both of $|\hat f(c)|$ and $|\hat f(c')|$ will
be smaller than $|\lambda_k|$.  From this observation, and a little work, 
we deduce that
\begin{eqnarray*}
\SUM_{a \in X} |\hat g(a)|^2\ &\leq&\ N^{-6} + 2k^{3\eps} k^{-(4+10\eps)(d-2)} \theta  \\
&&\ \ \ \ 
\times  \SUM_{-N/2 \leq b \leq N/2} \min(|I|,\ 
2^{-L} |I|^{-L+1} ||a/M - b/N||^{-L})^2.
\end{eqnarray*}
The factor $2k^{3\eps}$ is to account for the fact that given $c$ there are at most 
this many choices for $c'$ such that (\ref{dlc}) holds.
It is not difficult to see now that
\begin{equation} \label{axsum}
\SUM_{a \in X} |\hat g(a)|^2\ \leq\ 100 k^{-4(d-2) - \eps(10d-22)} \theta N^2. 
\end{equation}

\subsection{Selecting the right value for $m_1$, and therefore $M = m_1m_2$}

The value of $M$ that we will use will should be odd, coprime to $a_1,...,a_d$,
should satisfy
$$
k^\eps N\ <\ M\ \leq\ 4 k^\eps N,
$$
and should be factorable as
$$
M\ =\ m_1 m_2,\ {\rm where\ }{\rm gcd}(m_1,m_2)\ =\ 1,
$$
where $m_2$ is as we found previously.  Furthermore, $M$ should satisfy 
one more property, given as follows:
\bigskip

\noindent {\bf Correspondance property.}  
We want that $M$ satisfies the 
{\it correspondance property}, which is that for every pair of numbers
$$
x,\ y\ \in\ X^c
$$
and for every pair
$$
a_i,\ a_j\ \in\ \{a_1,...,a_d\},
$$
we have that 
$$
a_i x\ \equiv\ a_j y \pmod{M}\ \ \iff\ \ a_i x\ \equiv\ a_j y \pmod{m_2}.
$$
\bigskip

Another way of thinking of this property is as follows:  First, from the fact that 
$(m_1,m_2) = 1$, we may write 
$$
\Z_M\ =\ V\ +\ W,
$$
where $V$ and $W$ are subgroups given by
$$
V\ :=\ \{ m_1 x\ :\ 0 \leq x \leq m_2-1\},\ W\ :=\ \{ m_2 x\ :\ 0 \leq x \leq m_1-1\}.
$$
By the Chinese Remainder Theorem, every $a \in \Z_M$ may be written uniquely as
$$
a\ =\ v(a)\ +\ w(a),\ v(a) \in V,\ w(a) \in W.
$$
The correspondance property is then equivalent to saying that for 
$x,y \in X^c$,
\begin{equation} \label{correspondance_property}
a_i x\ \equiv\ a_i y \pmod{M}\ \ \iff\ \ a_i v(x)\ \equiv\ a_j v(y) \pmod{M}.
\end{equation}
Note that by the linearity of the projection maps $v$ and $w$ we have 
$v(a_ix) \equiv a_i v(x) \pmod{M}$.

\subsubsection{Proof that such $M$ exists}

As we have already selected $m_2$ in a previous subsection, it remains to 
find $m_1$.  To this end, fix a quadruple
$$
(b_i, b_j, a_u, a_v)\ \in\ \{b_1,...,b_k\}^2 \times \{a_1,...,a_d\}^2.
$$ 
Now suppose that $m_1$ is any integer satisfying
\begin{equation} \label{m1_range}
k^{-2-\eps} N\ <\ m_1 \ \leq\ 2k^{-2-\eps} N.
\end{equation}
We say that $m_1$ is ``good'' for this quadruple if for
every $x,y \in \Z_M$ satisfying
\begin{equation} \label{xy_close}
\left | \left | x/M - b_i/N \right | \right |,\ 
\left | \left | y/M - b_j/N \right | \right |\ \leq\ k^{3\eps}/M, 
\end{equation}
we have that 
$$
a_u x\ \equiv\ a_v y \pmod{m_2}\ \ \Longrightarrow\ \ a_u x\ \equiv\ a_v y \pmod{M}.
$$
Note that the reverse implication holds automatically.  

Clearly, if $m_1$ is ``good'' for every such quadruple $(b_i, b_j, a_u, a_v)$, then we
will have that $M=m_1m_2$ satisfies the correspondance property.
\bigskip

We now show that there can be few such integers $m_1$ that are ``bad'' for each
quadruple:  Suppose $m_1$ satisfies (\ref{m1_range}).    We will show that if $m_1$
satisfies the {\it separation property}, given in (\ref{separation_eitheror}),
then it must be ``good'' (in the sense above) for every quadruple. 

First, let us suppose that $x$ and $y$ satisfy (\ref{xy_close}).  Then,
$$
\left | \left |x/m_2\ -\ b_i m_1/N\right | \right |\ <\ k^{3\eps} /m_2,
$$
and the analogous inequality holds for $y/m_2$.  It is not difficult to see, then, that
$$
\left | \left | (a_u x - a_v y)/ m_2 \right | \right |\ =\ 
\left | \left | m_1(a_u b_i - a_v b_j)/N \right | \right |\ +\ \delta,
$$
where
$$
|\delta|\ \leq\ 2D k^{3\eps}/m_2. 
$$
Now, if we also add in the assumption that 
$$
a_ux\ \equiv\ a_v y \pmod{m_2},
$$
then we deduce that
$$
\left | \left | m_1(a_u b_i - a_v b_j)/N  \right | \right |\ \leq\ |\delta|.
$$
So, if $m_1$ is one of the integers satisfying (\ref{separation_eitheror}), then we 
are forced to have that
$$
a_u b_i\ \equiv\ a_v b_j \pmod{N}.
$$
But this, along with (\ref{xy_close}), implies that 
$$
\left | \left | (a_u x - a_v y)/ M \right | \right |\ \leq\ 2D k^{3\eps}/M.
$$
Since $a_u x - a_v y$ must be divisible by 
$$
m_2\ >\ 2 D k^{3\eps},
$$
we are forced to have 
$$
a_u x\ \equiv\ a_v y \pmod{M}.
$$

So, we have shown that if $m_1$ is one of the integers satisfying (\ref{m1_range}) 
and satisfying (\ref{separation_eitheror}) for every quadruple $(b_i,b_j,a_u,a_v)$,
then we will have that the {\it correspondance property} holds for $M=m_1m_2$.
Since we proved earlier that all but a fraction at most $k^{-2}$ of the
integers in (\ref{m1_range}) satisfy ({\ref{separation_eitheror}) for all
quadruples $(b_i,b_j,a_u,a_v)$, and since there are clearly fewer than this 
many integers $m_1$ having a common factor with 
${\rm lcm}(a_1,...,a_d,m_2)$, we must have that 
there exists $m_1$ satisfying all the above-mentioned properties 
(the correspondance property, as well as the coprimality conditions).

\subsubsection{A closer look at the correspondance property}
\label{closer_look_section}

In this subsection we use the fact that $M$ satisfies the correspondance property
to make two further deductions, listed below.
\bigskip

First, we split the subgroup $V$ into the two sets $V_1$ and $V_2$, where
$$
V_1\ :=\ \{v(a)\ :\ a \in X^c\},\ {\rm and\ } V_2\ :=\ V \setminus V_1.
$$
\bigskip

\noindent The first deduction is given as follows.
\begin{equation} \label{unique}
v \in V_1\ \ \Longrightarrow\ \  {\rm there\ exists\ unique\ }a \in X^c\ {\rm with\ } 
v(a)\ =\ v.
\end{equation}
To see that this holds, note that if there were two values $a,b \in X^c$ such that 
$$
v(a)\ \equiv\ v\ \equiv\ v(b) \pmod{M},
$$
then, say, 
$$
a_1 v(a)\ \equiv\ a_1 v(b) \pmod{M}\ \ \Longrightarrow\ \ 
a_1 a\ \equiv\ a_1 b \pmod{M},
$$
which then forces $a$ to equal $b$ mod $M$, hence $a$ is unique as claimed. 
\bigskip

The second deduction is that if 
$$
(a_1 v, ..., a_d v)\ \in\ V_1^d,
$$
then the unique $d$-tuple 
$$
(a_1 x_1,...,a_d x_d)\ \in\ (X^c)^d
$$
that maps as follows
$$
v\ :\ (a_1 x_1,...,a_d x_d)\ \to\ (v(a_1x_1),...,v(a_dx_d)) = (a_1 v,...,a_d v),
$$
must have the form
\begin{equation} \label{the_form}
(a_1 x_1,...,a_d x_d)\ =\ (a_1 b,...,a_d b),\ {\rm where\ } v = v(b).
\end{equation}
\bigskip

To make this deduction, we first observe that
$$
a_i v(a_1 x_1)\ \equiv\ a_i a_1 v\ \equiv\ a_1 a_i v\ \equiv\ a_1 v(a_i x_i) \pmod{M}.
$$
So, by (\ref{correspondance_property}), we deduce that
$$
x_1\ \equiv\ x_i \pmod{M}.
$$
But this means that the $d$-tuple $(a_1 x_1,...,a_d x_d)$ has the form 
(\ref{the_form}), as claimed.

\subsection{The auxilliary function $h$}

Using the function $g$ we now construct the auxiliary function $h$ as follows:
For a given translate $u \in \Z_M$, define
$$
V_u\ :=\ u + V\ =\ \{ u + v\ :\ v \in V\}.
$$
As usual, we associate to $V_u$ and $W$ the indicator functions $V_u(n)$ and 
$W(n)$ in the obvious way.  

Our function $h$ is to be defined as follows
$$
h\ :=\ h_u\ =\ (V_u g)*W.
$$
So,
$$
h(n)\ =\ \SUM_{a + b \equiv n \pmod{M}} V_u(a) g(a) W(b).
$$

We will need a formula for the Fourier transform of $h$:  We have that
for $a \notin V$, $\hat h(a) = 0$; on the other hand, if $a \in V$, then
\begin{eqnarray} \label{h_formula}
\hat h(a)\ &=&\  \widehat{(V_u g)}(a) \hat W(a) \nonumber \\
&=&\ m_2^{-1} \SUM_{x+y \equiv a \pmod{M}} \hat V_u(x) \hat g(y) \nonumber \\
&=&\ \SUM_{x \in W} e^{-2\pi i x u/M} \hat g(a+x).
\end{eqnarray}

\subsection{For some $u$, the function $h$ well approximates $g$
(in some sense)}

We claim that for each 
$$
a\ \in\ X^c 
$$
we will have that 
\begin{equation} \label{gh1}
\hat g(a)\ \ {\rm is\ very\ close\ to\ \ } e^{2\pi i w(a) u/M} \hat h(v(a)).
\end{equation} 
And, for $v \in V_2$ we will have 
\begin{equation} \label{gh2}
|\hat h(v)|\ \ {\rm is\ ``very\ small"}.
\end{equation}

To prove that there exists $u \in \Z_M$ so that 
both (\ref{gh1}) and (\ref{gh2}) hold, we will give an upper bound
for the following sum
\begin{equation} \label{sigma_def}
\Sigma\ =\ \SUM_{u \in \Z_M} 
\left ( \SUM_{a \in X^c} |\hat g(a)\ -\ e^{2\pi i w(a) u/M} \hat h(v(a))|^2
\ +\ \SUM_{v \in V_2} |\hat h(v)|^2\right ).
\end{equation}

To evaluate $\Sigma$ we first note from (\ref{h_formula}) that the first
term of (\ref{sigma_def}) equals
\begin{eqnarray} \label{uax1}
&& \SUM_{u \in \Z_M} \SUM_{a \in X^c} \left | \SUM_{x \in W \atop 
x \neq w(a)} e^{-2\pi i x u/M} \hat g(v(a) + x) \right |^2 \nonumber \\
&&\ \ \ \ =\  \SUM_{a \in X^c} \SUM_{x_1,x_2 \in W \atop 
x_1,x_2 \neq w(a)} \hat g(v(a) + x_1) \overline{\hat g(v(a) + x_2)}
\SUM_{u \in \Z_M} e^{-2\pi i u(x_1 - x_2)/M} \nonumber \\
&&\ \ \ \ =\  M \SUM_{a \in X^c} \SUM_{x \in W \atop x \neq w(a)}
|\hat g(v(a) + x)|^2.
\end{eqnarray}
The contribution of the second term in (\ref{sigma_def}) equals
\begin{eqnarray*}
\SUM_u \SUM_{v \in V_2} |\hat h(v)|^2\ &=&\ \SUM_u \SUM_{v \in V_2} 
\left | \SUM_{x \in W} e^{-2\pi i x u/M} \hat g(v + x) \right |^2 \\
&=&\ \SUM_{v \in V_2} \SUM_{x_1,x_2 \in W} \hat g(v+x_1)\overline{\hat g(v+x_2)} 
\SUM_u e^{-2\pi i u(x_1 - x_2)/M} \\
&=&\ M \SUM_{v \in V_2} \SUM_{x \in W} |\hat g(v+x)|^2.
\end{eqnarray*}
Combining this with (\ref{uax1}), the fact (\ref{unique}), and (\ref{axsum}), 
we find that 
$$
\Sigma\ =\ M \SUM_{a \in X} |\hat g(a)|^2\ <\ 100  k^{-4(d-2) - \eps(10d-22)} \theta M N^2.
$$
It follows that there exists
$u \in \Z_M$ such that 
\begin{eqnarray} \label{sigmafruit}
\SUM_{a \in X^c} |\hat g(a) - e^{2\pi i w(a) u/M} \hat h(v(a))|^2\ +\ 
\SUM_{v \in V_2} |\hat h(v)|^2\ <\ 100 k^{-4(d-2) - \eps(10d-22)} \theta N^2. \nonumber \\
\end{eqnarray}
We will use any such $u$ for our definition of $h$.
\bigskip

There are several conclusions that one can read off from this.  One such conclusion
is that
\begin{equation} \label{sigmafruit2}
a \in X^c\ \ \Longrightarrow\ \ |\hat g(a) - e^{2\pi i w(a) u/M} \hat h(v(a))|\ <\ 
10  k^{-2(d-2)-\eps(5d-11)} \theta^{1/2} N.
\end{equation}

\subsection{Relating the counting problem with $g$ to the counting problem
with $h$}

We now wish to consider the size of the error $E$ given by
\begin{eqnarray*}
E\ &=&\ \SUM_{x_1,...,x_d \in \Z_M \atop a_1 x_1 + \cdots + a_d x_d \equiv 0 \pmod{M}} 
g(x_1)\cdots g(x_d)\\
&&\hskip0.5in -\ 
\SUM_{x_1,...,x_d \in \Z_M \atop a_1 x_1 + \cdots + a_d x_d \equiv 0 \pmod{M}} 
h(x_1)\cdots h(x_d).
\end{eqnarray*}

This can be expressed in terms of Fourier coefficients as
\begin{eqnarray} \label{Eformula}
E\ &=&\ M^{-1} \SUM_b \hat g(a_1b )\cdots \hat g(a_db )\ -\ \hat h(a_1b)\cdots \hat h(a_d b) 
\nonumber \\
&=&\ M^{-1} \SUM_b \hat g(a_1 b) \cdots \hat g(a_d b)\ -\ M^{-1} 
\SUM_{v \in V} \hat h(a_1 v)\cdots \hat h(a_d v).
\end{eqnarray}

\subsubsection{Contribution of $b$ where $\hat g$ is small}

Let us now consider the contribution to this first sum, all those $b$ where
\begin{equation} \label{bsatisfy}
{\rm either\ }a_1 b \in X,\ {\rm or\ } a_2 b \in X,\ ...,\ {\rm or\ } a_d b \in X.
\end{equation}
Using (\ref{axsum}) the contribution of those $b$ with $a_1 b \in X$ 
can be bounded from above by
\begin{eqnarray} \label{error1}
&& M^{-1} \left ( \SUM_{x \in X} |\hat g(x)|^2 \right )^{1/2} 
\left ( \SUM_b |\hat g(b)|^{2d-2} \right )^{1/2} \nonumber \\
&&\hskip0.5in <\ 10 k^{-2(d-2)-\eps(5d-11)} \theta^{1/2} \left ( 
\SUM_b |\hat g(b)|^{2d-2} \right )^{1/2}
\end{eqnarray}
To handle this last factor, we first note that
$$
\hat g(b)^{d-1}\ =\ \SUM_{n \in \Z_M} \nu(n) e^{2\pi i bn/M} , 
$$
where
$$
\nu(n)\ =\ \SUM_{y_1 + \cdots + y_{d-1} \equiv n \pmod{M}} 
g(y_1)\cdots g(y_{d-1}).
$$
From the fact that 
$$
\hat g(0)\ \leq\ \hat f(0)\ =\ \theta N,
$$
one can show that the sum over $\nu(n)^2$ is maximized if $\nu(n)$ is supported
on an interval of size $\theta N$, and at each such place $\nu(n)$ has value at most 
$(\theta N)^{d-2}$.  So,
$$
\SUM_n \nu(n)^2\ \leq\ (\theta N)^{2d-3},
$$
and this, along with Parseval, implies that
$$
\SUM_b |\hat g(b)|^{2d-2}\ \leq\ M (\theta N)^{2d-3}\ \leq\ 4k^\eps \theta^{2d-3} N^{2d-2}.
$$
It follows that (\ref{error1}) is bounded from above by
$$
20 k^{-2(d-2) - \eps(5d-23/2)} (\theta N)^{d-1}.
$$
Recalling that this is only an upper bound for the contribution to the first 
expression in (\ref{Eformula}) with $a_1 b \in X$, we find that the total contribution
to this expression with $b$ satisfying (\ref{bsatisfy}) is bounded from above by
$$
20 d k^{-2(d-2)-\eps(5d-23/2)} (\theta N)^{d-1}.
$$

\subsubsection{Contribution of $b$ where $\hat h$ is small}

Now we consider the contribution to the second sum of (\ref{Eformula}) of
all those $v \in V$ such that 
\begin{equation} \label{allv}
{\rm either\ }a_1 v,\ {\rm or\ } a_2 v,\ ...,\ {\rm or\ } a_d v\ \in\ V_2.
\end{equation}
First, let us consider the contribution of those terms with $a_1 v \in V_2$:  Using
(\ref{sigmafruit}), this can be bounded from above by
\begin{eqnarray} \label{herror1}
&& M^{-1} \left ( \SUM_{v \in V_2} |\hat h(v)|^2 \right )^{1/2}
\left ( \SUM_{v \in V} |\hat h(v)|^{2d-2} \right )^{1/2} \nonumber \\
&& \hskip1in <\ 10 k^{-2(d-2) - \eps(5d-11)} \theta^{1/2}
\left ( \SUM_{v \in V} |\hat h(v)|^{2d-2} \right )^{1/2}.
\end{eqnarray}

To bound this remaining sum over $|\hat h(v)|^{2d-2}$ from above, we first observe
from (\ref{sigmafruit2}), along with the fact that $\hat g(0) \leq \hat f(0)$, that 
\begin{eqnarray} \label{hupperdensity}
\hat h(0)\ &\leq&\ \hat g(0)\ +\ |\hat h(0)\ -\ \hat g(0)| \nonumber \\
&\leq&\ \theta N\ +\ 10 k^{-2(d-2)-\eps(5d-11)} \theta^{1/2} N \nonumber \\
&<&\ \theta (1 + d^{-1}) N,
\end{eqnarray}
which follows since
$$
k\ >\ \theta^{-1/4(d-2)}. 
$$
Now write
$$
\hat h(v)^{d-1}\ =\ \SUM_{n \in \Z_M} \nu'(n) e^{2\pi i vn/M},
$$
From our upper bound on $\hat h(0)$, one can easily deduce that
$$
\SUM_{n \in \Z_M} \nu'(n)^2\ <\ 10 (\theta N)^{2d-3} 
$$
So, using Parseval this gives
$$
\SUM_{v \in V} |\hat h(v)|^{2d-2}\ \leq\ 10 M (\theta N)^{2d-3};
$$
and so, the quantity in (\ref{herror1}) is bounded from above by
$$
100 k^{-2(d-2) - \eps(5d-23/2)}  (\theta N)^{d-1}.
$$
This was just the contribution to (\ref{Eformula}) of all places $v \in V$ 
where $a_1 v \in V_2$.  The contribution of all $v$ satisfying (\ref{allv}) is 
thus 
$$
100 d k^{-2(d-2)-\eps(5d-23/2)} (\theta N)^{d-1}.
$$

\subsubsection{Comparison of the main terms}

Now we are left to consider the contribution to (\ref{Eformula}) of all those
$b \in \Z_M$ and $v \in V$ satisfying 
$$
(a_1 b, a_2 b, ..., a_d b)\ \in\ (X^c)^d,\ {\rm and\ } (a_1 v, a_2 v, ..., a_d v)\ \in\ V_1^d.
$$

We note that from our deduction in sub-subsection \ref{closer_look_section}
that there is a one-to-one correspondance between the set of such $b$ and the 
set of such $v$ making so that these $d$-tuples are in $(X^c)^d$ and
$V_1^d$, respectively.  Indeed, the correspondance is such that 
$v = v(b)$.  So, we may index both types of $d$-tuples by certain elements $b \in \Z_M$.
\bigskip

For each of these $b$, let us now consider
\begin{equation} \label{two_terms}
\hat g(a_1 b)  \hat g(a_2 b) \cdots \hat g(a_d b)\  -\  
\hat h(a_1 v(b)) \hat h(a_2 v(b)) \cdots \hat h(a_d v(b)).
\end{equation}

Applying (\ref{sigmafruit2}) we can replace each of the factors $\hat g(a_i b)$
of (\ref{two_terms}) with $e^{2\pi i w(a_i b) u/M} \hat h(a_i v(b))$ plus a 
small error, say
\begin{eqnarray}
\hat g(a_1 b) \cdots \hat g(a_d b)\ &=&\ \prod_{i=1}^d \biggl (
e^{2\pi i w(a_i b) u/M} \hat h(a_i v(b)) + F_i \biggr ) \\
&=&\ \prod_{i=1}^d \biggl ( \hat h(a_i v(b)) + F_i' \biggr ),
\end{eqnarray}
where the $F_i$ and $F_i'$ are the errors, with $|F_i| = |F_i'|$.  
Note that the reason we can get rid of the roots of unity factors is that 
$$
a_1 + \cdots + a_d\ =\ 0\ \ \Longrightarrow\ \ w(a_1 b) + \cdots + w(a_d b)\ \equiv\ 0
\pmod{M},
$$
since $w$ is a linear map.

If we expand out this product we get
$$
\hat h(a_1 v(b)) \cdots \hat h(a_d v(b))
$$
as a ``main term'', and then we get $2^d - 1$ ``error terms'' in all.  Suppose we 
fix one of these error terms, and sum its absolute value over all $b \in \Z_M$.  
It is easy to see, using
(\ref{hupperdensity}), (\ref{sigmafruit2}), Parseval's identity, and Holder's inequality, 
that an upper bound for the result is
\begin{eqnarray*}
\sup_i |F_i'| \SUM_b |\hat h(b)|^{d-1}\ &\leq&\ (\sup_i |F_i'|) |\hat h(0)|^{d-3} 
\SUM_b |\hat h(b)|^2 \\ 
&\leq&\ 30 k^{-2(d-2)-\eps(5d-11)} \theta^{d-3/2} N^{d-1} M. 
\end{eqnarray*}

That was the error arising from just one of the $2^d - 1$ terms, and before multplying
through by the $M^{-1}$ out front of (\ref{Eformula}).  In all, then, the
contribution of these ``main term'' errors to (\ref{Eformula}) is at most 
$$
2^{d+5} k^{-2(d-2)-\eps(5d-11)} \theta^{-1/2} (\theta N)^{d-1}.
$$

\subsubsection{The main term for $h$, and the conclusion of the proof}

To finish the proof of our theorem, we will bound the following from below
\begin{equation} \label{h_count_end}
M^{-1} \SUM_{v \in V} \hat h(a_1 v) \cdots \hat h(a_d v)\ =\ 
\SUM_{x_1,...,x_d \in \Z_M \atop a_1 x_1 + \cdots + a_d x_d 
\equiv 0 \pmod{M}} h(x_1)\cdots h(x_d).
\end{equation}

To do this we will use the fact that $h$ is translation-invariant by elements 
of $w \in W$, along with a lower bound for $\hat h(0)$.   This lower bound is
a companion to (\ref{hupperdensity}), and is proved in exactly the same way:
From (\ref{gmass}), (\ref{sigmafruit2}), and the fact that
$$
k\ >\ \theta^{-1/4(d-2)},
$$
we deduce that
\begin{eqnarray*}
\hat h(0)\ \geq\ \hat g(0) - |\hat g(0) - \hat h(0)|\ &\geq&\ k^{-\eps}(1 - 1/d) \theta N \\
&\geq&\ 4^{-1} k^{-2\eps} (1 - 1/d) \theta M.
\end{eqnarray*}
Another fact we will use is that $h$ is translation-invariant by elements 
of $W$; that is, we will use the fact that for $w \in W$,
$$
h(n+w)\ =\ [(V_u g)*W](n+w)\ =\ [(V_u g)*W](n)\ =\ h(n).
$$ 
The way that this helps us to bound (\ref{h_count_end}) from below is that for 
each $n \in \Z_M$ and for each $d-1$ tuple 
$$
w_1,\ ...,\ w_{d-1}\ \in\ W,
$$
if we let 
$$
w_d\ \equiv\ - a_d^{-1} (a_1 w_1 + \cdots + a_{d-1} w_{d-1}) \pmod{M},
$$
then we have that 
$$
w_d\ \in\ W,\ {\rm and\ }
a_1 w_1 + \cdots + a_d w_d\ \equiv\ 0 \pmod{M}.
$$

This implies in particular that for all $n \in \Z_M$,
$$
a_1 (n + w_1) + \cdots + a_d (n + w_d)\ \equiv\ 0 \pmod{M},
$$
and 
$$
h(n)^d\ =\ h(n+w_1) \cdots h(n+w_d).
$$
So, a lower bound for (\ref{h_count_end}) is given by
\begin{eqnarray} \label{last_line}
|W|^{d-1} \SUM_{n \in V} h(n)^d\ &=&\ |W|^{d-2} \SUM_{n \in \Z_M} h(n)^d \nonumber \\
&\geq&\ 4^{-d} k^{-d\eps} (1 - 1/d)^d  \theta^d M |W|^{d-2}  \nonumber \\
&\geq&\ e^{-1} 4^{-d} k^{-d\eps} \theta^d M^{d-1} m_2^{-d+2} \nonumber \\
&\geq&\ e^{-1} 4^{-d} k^{-2(d-2) - 2\eps(d-3/2)} \theta (\theta N)^{d-1} 
\end{eqnarray}

Now we compare this with the sum of all the ``error terms'' accumulated from
the previous subsubsections.  Since 
$$
k\ >\ 1000^{d\eps^{-1}} \theta^{-(\eps d)^{-1}},\ {\rm and\ } d\ \geq\ 3,
$$ 
these all sum to at most 
\begin{eqnarray*}
&& 2^{d+5} k^{-2(d-2)-\eps(5d-23/2)} \theta^{-1/2} (\theta N)^{d-1} \\
&& \hskip1in +\ 100 d k^{-2(d-2)-\eps(5d-23/2)} (\theta N)^{d-1} \\
&& \hskip1in +\ 20 d k^{-2(d-2)-\eps(5d-23/2)} (\theta N)^{d-1} \\
&&\ \leq\ 10^{-1} 4^{-d} k^{-2(d-2) - 2\eps(d-3/2)} \theta (\theta N)^{d-1}.
\end{eqnarray*}
This last quantity is at most half that in the last line of (\ref{last_line}); and 
so, 
$$
\SUM_{x_1,...,x_d \in \Z_N \atop a_1 x_1 + \cdots + a_d x_d \equiv 0 \pmod{N}}
f(x_1)\cdots f(x_d)\ >\ 10^{-1} 4^{-d} k^{-2(d-2) - 2\eps(d-3/2)} \theta (\theta N)^{d-1}, 
$$
as claimed.

\end{document}